\documentclass{article}      

\usepackage{epsfig}

\title{NON-COMMUTATIVE $q$-BINOMIAL  FORMULA}

\author{Sengul Nalci and Oktay K. Pashaev \\Department of Mathematics, Izmir Institute of Technology \\ Urla-Izmir, 35430, Turkey}

\begin{document}
\newcommand{\be}{\begin{equation}}
\newcommand{\ee}{\end{equation}}
\newcommand{\bea}{\begin{eqnarray}}
\newcommand{\eea}{\end{eqnarray}}
\newcommand{\disp}{\displaystyle}
\newcommand{\la}{\langle}
\newcommand{\ra}{\rangle}

\newtheorem{thm}{Theorem}[subsection]
\newtheorem{cor}[thm]{Corollary}
\newtheorem{lem}[thm]{Lemma}
\newtheorem{prop}[thm]{Proposition}
\newtheorem{defn}[thm]{Definition}
\newtheorem{rem}[thm]{Remark}
\newtheorem{prf}[thm]{Proof}

\maketitle


\begin{abstract}
In this paper, we found new $q$-binomial formula for $Q$-commutative operators. Expansion coefficients in this formula are given by $q$-binomial coefficients with two bases $(q,Q),$ determined by $Q$-commutative $q$-Pascal triangle. Our formula generalizes all well-known binomial formulas in the form of Newton, Gauss, symmetrical, non-commutative and Binet-Fibonacci binomials. By our non-commutative $q$-binomial, we introduce  $q$-analogue of function of two non-commutative  variables, which could be used in study of non-commutative $q$-analytic functions and non- commutative $q$-traveling waves.
\end{abstract}

\section{Introduction}

The Newton's Binomial Formula for positive integer $n$ is given in the following form
\bea (x+y)^n = \sum_{k=0}^n {n \choose k } x^{n-k} y^k, \eea
where \bea {n \choose k }= \frac{n!}{k! (n-k)!}\nonumber \eea denotes the binomial coefficients.

The $q$-analogue of the Newton's Binomial formula is given for the $q$-analogue of binomial $(x-a)^n,$ called the $q$-binomial.

\begin{defn} The $q$-analogue of $(x-a)^n$ is the polynomial
{\small \[ (x-a)_{q}^n = \left \{ \begin{array} {ll} 1 & \mbox{if $n=0$}, \\
 (x-a) (x-q a) (x-q^2 a)...(x-q^{n-1} a)
        & \mbox {if $n\geq 1$ } \end{array}
        \right.\] }
\end{defn}
For this $q$-binomial the
Gauss's Binomial formula for \textbf{commutative} $x$ and $a$ $(xa=ax)$ is written in the following form \cite{Kac et al.}
\bea (x+a)_q^n = \sum_{k=0}^n {n \brack k }_q q^{\frac{k(k-1)}{2}}x^{n-k} a^k,  \label{gaussbinom}\eea
where
the $q$-binomial coefficients are $${n \brack k }_q= \frac{[n]_q!}{ [k]_q! [n-k]_q!},$$ and $$[n]_q!=[1]_q  [2]_q ...  [n]_q\,, \,\,\,\,\,\,\,\,\,\,\,\,\,\,\, [n]_q= \frac{q^n-1}{q-1}.$$

In addition to the Gauss's Binomial formula, the \textbf{non-commutative } Binomial formula
for the \textbf{ $q$-commutative} x and y $(y x=q x y)$ is valid \cite{Kac et al.}
\bea (x+y)^n =\sum_{k=0}^n {n \brack k }_q x^k y^{n-k}, \eea
 where $q$ is a number, commutating with $x$ and $y$; $x q=q x$ and $y q= q y.$

In the present paper, we derive new type of non-commutative $q$-binomial for $q$-binomial with $q$-commutative entries. We show that this type of formula is naturally related with $q$-calculus with two bases $(Q,q).$

\section{Q-commutative q-Binomial Formula}
In several physical and mathematical problems we encounter to find a  new type of  binomial formula.
The $q$-analogue of some known PDE like $q$-wave and $q$-Burgers equation include $Q$-derivative operator $D_Q$ and the dilatation operator $M_Q$ which are non-commutative, but $Q$-commutative operators. Solutions of these equations are represented as $q$-binomials in terms of non-commutative operators \cite{Nalci et al.}. Non-commutative versions of $q$-traveling waves \cite{cNalci et al.} with non-commutative $x$ and $t$ provides non-commutative space-time equations. From another side, extension of quantum groups to two basis $(Q,q)$ require $q$-calculus with multiple $q$-numbers. Quantum $q$-oscillator with such symmetries has been discussed in \cite{Pashaev and Nalci}. It turns out that non-commutative $q$-binomial formulas are naturally described in terms of such two base calculus.

Firstly, we note that in the standard notation of $q$-binomial (we using notation from Kac \cite{Kac et al.})
\bea (x+y)_q^n= (x+y) (x+q y) (x+q^2 y)...(x+q^{n-1} y), \,\,\,\,\, n=1,2,..\eea
 applied to the noncommutative operators $x$ and $y$,  we should distinguish the direction of multiplication. So  we introduce the following notation for two different cases of order \cite{Nalci et al.}
\bea (x+y)^n_{<q}\equiv (x+y) (x+q y) (x+q^2 y)...(x+q^{n-1}y ) \label{kb}\eea and
\bea (x+y)^n_{>q} \equiv (x+q^{n-1}y)... (x+qy) (x+y).\label{bk} \eea
Before introducing the $q$- Binomial formula for  $Q$-commutative operators, we briefly give the definition of multiple $q$-numbers :

Multiple $q$-Number is defined by a basis vector $\overrightarrow{q}$ with coordinates $q_1,q_2,...,q_N$ as a matrix $q$-number,
\bea [n]_{{q_i},{q_j}} \equiv \frac{q_i^n-q_j^n}{q_i-q_j}=[n]_{{q_j},{q_i}} , \eea which is symmetric.

\begin{thm}
Let  $x$ and $y$ are  $Q$-commutative operators, $y x=Q x y,$ then the $q$-Binomial formula is valid
\bea (x+y)^n_{<q}= \sum_{k=0}^n {n \brack k }_{Q,q} q^{\frac{k(k-1)}{2}} x^{n-k} y^k,\eea where 
$(Q,q)$- binomial coefficients are defined as \bea {n \brack k}_{Q,q}= \frac{[n]_{Q,q}!}{[n-k]_{Q,q}![k]_{Q,q}!},\label{multiplecombinatorix}\eea and the $(Q,q)$-numbers are \bea [n]_{Q,q}= \frac{Q^n-q^n}{Q-q}. \label{multipleqnumber}\eea
\end{thm}

\begin{prf}
Here we give two differen proofs of this theorem. The first one reduces the problem to solution of system of difference equations. The second one is based on the method of mathematical induction. 

To find expansion of  $q$-polynomials in terms of $x$ and $y$ powers we 
suppose the following expansion
 \bea (x+y)^n_{< q}=\sum_{k=0}^n {n \brace k}_{Q,q} x^{n-k} y^k, \label{noncomgaus}\eea
 where  \bea (x+y)^n_{< q}= (x+y) (x+ q y) (x+q^2 y)... (x+q^{n-1} y),\nonumber \eea
and  ${n \brace k}_{Q,q}$- denote unknown coefficients, depending on $k,n,q$ and $Q.$
Then by definition ,
\bea (x+y)^{n+1}_{<q}= (x+y)^n _{<q}  (x+q^n y).\nonumber \eea
Expanding both sides
\bea  \sum_{k=0}^{n+1} {{n+1} \brace k}_{Q,q} x^{n-k+1} y^k&=& \sum_{k=0}^n {n \brace k}_{Q,q} x^{n-k} y^k (x+q^n y)\nonumber \\
&=& \sum_{k=0}^n {n \brace k}_{Q,q} x^{n-k} y^k x + \sum_{k=0}^n {n \brace k}_{Q,q} q^n x^{n-k} y^{k+1} \nonumber \\
&=& \sum_{k=0}^n {n \brace k}_{Q,q} Q^k x^{n-k+1} y^k  + \sum_{k=0}^n {n \brace k}_{Q,q} q^n x^{n-k} y^{k+1} \nonumber \\
&=& \sum_{k=0}^n {n \brace k}_{Q,q} Q^k x^{n-k+1} y^k + \sum_{k=1}^{n+1} {n \brace {k-1}}_{Q,q} q^n x^{n-k+1} y^k  \nonumber \eea

from the above equality we have the following recursion formulas :
\bea k= 0 \Rightarrow \,\,\,\,\,\, {{n+1} \brace 0 }_{Q,q}&=&{n \brace 0 }_{Q,q}, \nonumber \\
k= n+1 \Rightarrow \,\,\,\,\,\, {{n+1} \brace {n+1}}_{Q,q}&=&q^n {n \brace n }_{Q,q}, \nonumber \\
1\leq k \leq n \Rightarrow \,\,\,\,\,\, {{n+1} \brace k }_{Q,q}&=&Q^k {n \brace k}_{Q,q}+ q^n {n \brace {k-1} }_{Q,q}, \label{recursion}\eea
where we choose \bea {n \brace b }_{Q,q}=0 \,\,\, \rm{if} \,\,\, b<0 \,\,\, \rm{and} \,\,\, b>n. \nonumber \eea
Suppose the unknown binomial coefficient factor ${n \brace k }_{Q,q}$ can be written in terms of the known combinatorial coefficient ${n \brack k }_{Q,q}$ with multiplication factor as
\bea {n \brace k }_{Q,q}= q ^{t(n,k)}{n \brack k }_{Q,q}, \eea
where ${n \brack k }_{Q,q}$ is $(q,Q)$- combinatorial coefficient (\ref{multiplecombinatorix}).

Substituting this relation to  (\ref{recursion}) and using the following relation for $q$-multiple binomial coefficients, by choosing $q_i=Q \,\,\,\, q_j=q,$
\bea { n \brack k}_{q_i,q_j}&=& \frac{q_j^k [n-1]_{q_i,q_j}!}{[k]_{q_i,q_j}! [n-k-1]_{q_i,q_j}!} + \frac{q_i^{n-k} [n-1]_{q_i,q_j}!}{[n-k]_{q_i,q_j}![k-1]_{q_i,q_j}!} \nonumber \\
&=& q_j^k { n-1 \brack k}_{q_i,q_j} + q_i^{n-k} { n-1 \brack k-1}_{q_i,q_j} ,\label{qpascal1}\\
&=& q_i^k { n-1 \brack k}_{q_i,q_j} + q_j^{n-k} { n-1 \brack k-1}_{q_i,q_j},\label{qpascal2} \eea
we have following expression
\bea Q^k q^{t(n+1,k)}{n \brack k }_{Q,q} + q^{n+1-k+t(n+1,k)}{n \brack {k-1}}_{Q,q}&=&Q^k q^{t(n,k)}{n \brack k }_{Q,q}\nonumber \\
&+&q^{n+t(n,k-1)}{n \brack {k-1} }_{Q,q}  \nonumber \eea
By equating terms with the same power of $q$ and $Q,$ we obtain two difference equations
\bea t(n+1,k)&=& t(n,k) \nonumber \\
t(n,k)&=& t(n,k-1)+k-1 \eea with the initial conditions \bea t(0,0)= t(1,0)=t(1,1)=0.\eea
From the first equation for $k=0$ we have $t(n+1,0)=t(n,0).$
So, if $n=1\Rightarrow \,\, t(2,0)= t(1,0)=0,$ which means that $t(n,0)=0.$
By using the second equation we easily write
\bea t(n,1)&=&t(n,0)=0 \nonumber \\
 t(n,2)&=& t(n,1)+1=1 \nonumber \eea
 \bea
t(n,3)&=& t(n,2)+2=1+2 \nonumber \\
t(n,4)&=& t(n,3)+2=1+2+3 \nonumber \\
... \nonumber \\
t(n,k)&=& 1+2+3+...+(k-1)= \frac{k(k-1)}{2}. \nonumber \eea
Therefore, the solution of the above system is
\be t(n,k)= \frac{k(k-1)}{2}.\ee
Hence, we obtain the $q$- Binomial formula for  $Q$-commutative $x$ and $y$  in the form
\bea (x+y)^n_{<q}= \sum_{k=0}^n {n \brack k }_{Q,q} q^{\frac{k(k-1)}{2}} x^{n-k} y^k,\eea where $y x=Q x y,$ and
\bea {n \brack k}_{Q,q}= \frac{[n]_{Q,q}!}{[n-k]_{Q,q}![k]_{Q,q}!}, \,\,\,\,\,\, [n]_{Q,q}= \frac{Q^n-q^n}{Q-q}.\nonumber \eea
\end{prf} \hfill \rule{1.6ex}{1.6ex}

\begin{figure}[h]
\begin{center}
{\includegraphics[width=3.5in]{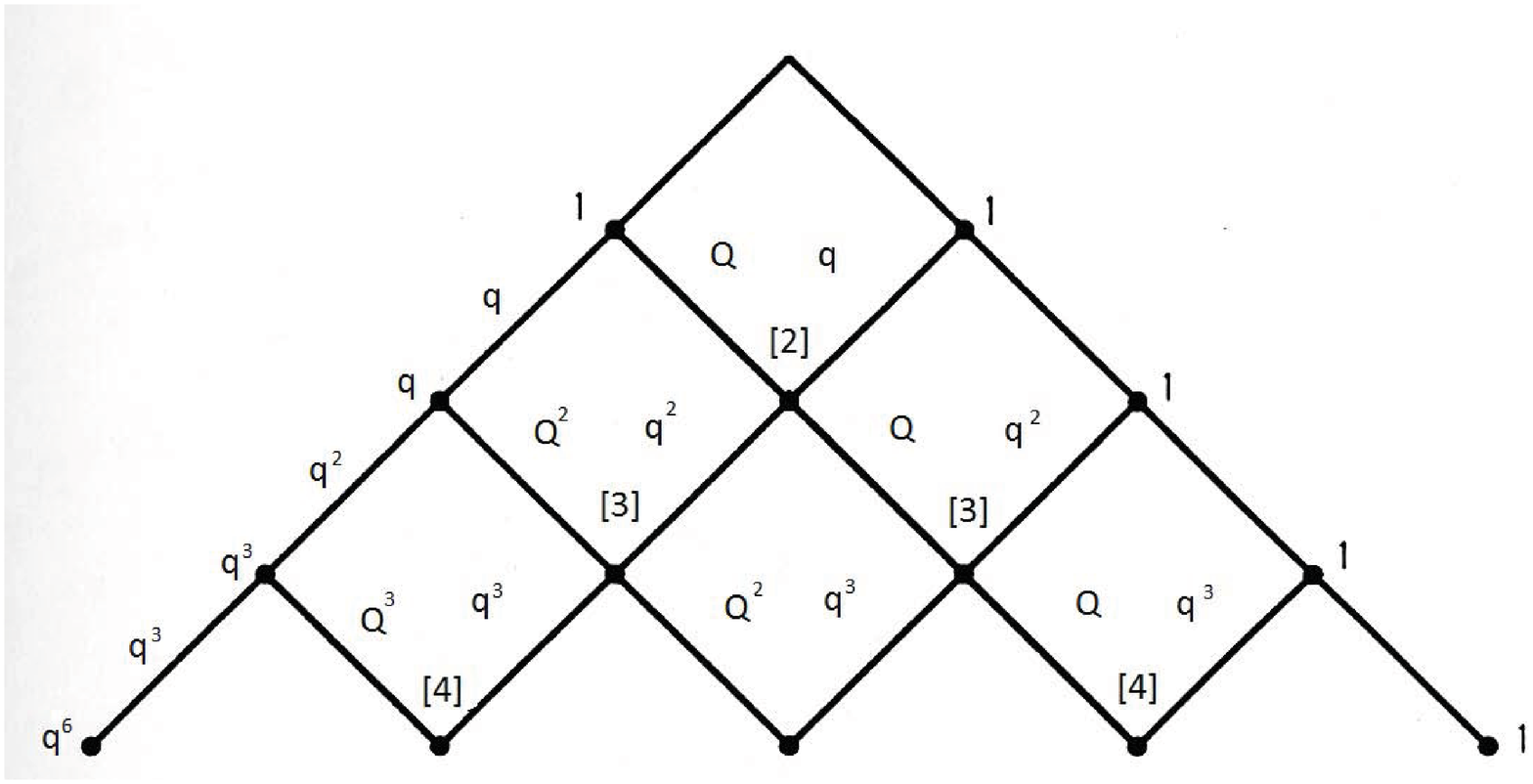}\\}
\caption{Q-commutative q-Pascal triangle}\label{qqpascal}
\end{center}
\end{figure}

\begin{prf}
It is instructive now to prove this Binomial formula by using mathematical induction.
We have
\bea (x+y)^{n+1}_{<q} &=& (x+y)^n_{<q} (x+q^n y) \nonumber \\
&=& \sum_{k=0}^n {n \brack k }_{Q,q} q^{\frac{k(k-1)}{2}} x^{n-k} y^k (x+q^n y)\nonumber \\
&=&  \sum_{k=0}^n {n \brack k }_{Q,q} q^{\frac{k(k-1)}{2}} x^{n-k} y^k x + \sum_{k=0}^n {n \brack k }_{Q,q} q^{\frac{k(k-1)}{2}} q^n x^{n-k} y^{k+1}.\nonumber \eea
From the $Q-$ commutativity relation $y x=Q x y,$ we get $y^k x= Q^k x y^k$ and the above expression is written as follows
\bea (x+y)^{n+1}_q &=&  \sum_{k=0}^n {n \brack k }_{Q,q} q^{\frac{k(k-1)}{2}}Q^k x^{n-k+1} y^k  + \sum_{k=0}^n {n \brack k }_{Q,q} q^{\frac{k(k-1)}{2}} q^n x^{n-k} y^{k+1} \nonumber \\
&=&  \sum_{k=0}^n {n \brack k }_{Q,q} q^{\frac{k(k-1)}{2}}Q^k x^{n-k+1} y^k  + \sum_{k=1}^{n+1} {n \brack {k-1} }_{Q,q} q^{\frac{(k-1)(k-2)}{2}} q^n x^{n-k+1} y^{k} \nonumber \\
&=& {n \brack 0 }_{Q,q} x^{n+1} + {n \brack n }_{Q,q} q^{\frac{n(n+1)}{2}} y^{n+1}\nonumber \\
&+& \sum_{k=1}^{n} \left({n \brack k }_{Q,q} q^{\frac{k(k-1)}{2}}Q^k + {n \brack {k-1} }_{Q,q} q^{\frac{(k-1)(k-2)}{2}} q^n \right) x^{n-k+1} y^k. \label{prove}\eea
  By choosing $q_i=Q $ and $q_j=q$  for binomial coefficients (\ref{qpascal2}), we have the following $(Q,q)$-Pascal triangle relation (see Figure 1).
\bea Q^k {n \brack k }_{Q,q}={{n+1} \brack k }_{Q,q}- q^{n+1-k}{n \brack {k-1} }_{Q,q}.\nonumber \eea  By substituting  this relation into equation (\ref{prove}) we have desired result
\bea(x+y)^{n+1}_q &=& {n \brack 0 }_{Q,q} x^{n+1} + {n \brack n }_{Q,q} q^{\frac{n(n+1)}{2}} y^{n+1}+ \sum_{k=1}^{n} {n \brack k }_{Q,q} q^{\frac{k(k-1)}{2}}x^{n-k+1} y^k \nonumber \\
&=& \sum_{k=0}^{n+1}{{n+1} \brack k }_{Q,q} q^{\frac{k(k-1)}{2}} x^{n-k+1} y^k. \eea
\end{prf}
 \hfill \rule{1.6ex}{1.6ex}

\textbf{Example:}
We consider the $Q$- derivative operator $D_Q= \frac{M_Q -1}{x(Q-1)},$ where $M_Q=Q^{x \frac{d}{d x}}$ is the $Q$-dilatation operator. These operators are $Q$-commutative \bea D_Q M_Q=Q M_Q D_Q.\nonumber \eea
Then according to the theorem, we have the binomial expansion
\bea (M_Q+D_Q)^n_{<q} &=& (M_Q+D_Q)(M_Q+q D_Q)(M_Q+q^2 D_Q)...(M_Q+q^{n-1}D_Q) \nonumber \\
&=&  \sum_{k=0}^n {n \brack k }_{Q,q} q^{\frac{k(k-1)}{2}} M_Q^{n-k} D_Q^k .\eea

\textbf{Example:}
The unitary operator $D(\alpha)= e^{\alpha a^+ - \bar{\alpha}a}$ for the Heisenberg-Weyl group is the generating operator for Coherent states \bea | \alpha >= D(\alpha)|0\ra, \nonumber \eea
where $\alpha$ is complex parameter. These operators satisfy the relation \bea D(\alpha) D(\beta)= e^{2i \Im (\alpha \bar{\beta})} D(\beta)D(\alpha),\nonumber \eea hence they are $Q$-commutative with $Q=e^{2i \Im (\alpha \bar{\beta})}.$ Therefore, applying the above theorem we have the next operator $q$-Binomial expansion
\bea \left(  D(\beta)+ D(\alpha)\right)^n_{<q} = \sum_{k=0}^n {n \brack k}_{Q, q} q^{\frac{k(k-1)}{2}}D^{n-k}(\beta)D^k(\alpha).  \eea
If we apply this expansion to vacuum state $|0 \ra $, then we get expansion in superposition of Coherent states 
\bea \left(  D(\beta)+ D(\alpha)\right)^n_{<q}|0 \ra &=& \sum_{k=0}^n {n \brack k}_{Q, q} q^{\frac{k(k-1)}{2}}D^{n-k}(\beta) | k \alpha \ra \nonumber \\
&=& \sum_{k=0}^n {n \brack k}_{Q, q} q^{\frac{k(k-1)}{2}} e ^{i \delta} | (n-k) \beta +k \alpha \ra ,\eea
where \bea \delta= 2 \Im \left[(n-k)k \alpha \bar{\beta}\right].\nonumber \eea

In the 
above theorem we have derived $Q$-commutative (\ref{kb}) $q$-binomial formula for ordered product $<q.$ Now we like to construct similar formula for opposite ordered product $>q$ (\ref{bk}).
\begin{prop}
 The following relation for two different ordered products is valid
\bea \prod_{k=0}^N (x+q^k y )_{<q}=\prod_{k=0}^N (x+Q^{N-2k} q^k y )_{>q}, \label{relationbk}\eea
where $yx=Qxy$ and
\bea & &\prod_{k=0}^N (x+q^k y )_{<q}=(x+y) (x+ q y) (x+q^2 y)... (x+q^N y)= (x+y)^{N+1}_{<q}, \nonumber \\
& & \prod_{k=0}^N (x+Q^{N-2k} q^k y)_{>q}=(x+q^N y)(x+q^{N-1}y)...(x+q y) (x+y).\nonumber \eea
\end{prop}
\begin{prf}
This formula can be proved by the method of mathematical induction.
\bea N=1\Rightarrow \,\,\, \prod_{k=0}^1 (x+q^k y)_{<q}&=& (x+y) (x+qy)=x^2 +qQ^{-1} yx+Qxy+qy^2 \nonumber \\ &=& (x+Q^{-1}q y)(x+Q y)
=\prod_{k=0}^1 (x+Q^{1-2k} q^k y )_{>q} \nonumber \eea
and we suppose that the formula is true for some $N.$ Let us show that it is also valid for $N+1:$
\bea \prod_{k=0}^{N+1} (x+q^k y)_{<q} &=& \prod_{k=0}^N (x+q^k y)_{<q}(x+q^{N+1} y)\nonumber \\
&=& \prod_{k=0}^N (x+Q^{N-2k} q^k y)_{>q}(x+q^{N+1} y)\nonumber \\
&=& (x+Q^{-N} q^N y) ...(x+Q^{(N-2)} q y)(x+Q^{N} y)(x+q^{N+1} y).\nonumber\eea
By using  equality
\bea (x+q^m y) (x+q^k y)= (x+Q^{-1} q^k y)(x+Q q^m y),\nonumber \eea
we move the last term to the left-end by commutating with every term of the product
\bea \prod_{k=0}^{N+1} (x+q^k y)_{<q} &=& (x+Q^{-N} q^N y) ...(x+Q^{(N-2)} q y)(x+Q^{-1}q^{N+1} y)(x+Q^{N+1}y)\nonumber \\
&=& (x+Q^{-N} q^N y) ...(x+Q^{-2}q^{N+1} y)(x+Q^{(N-1)} q y)(x+Q^{N+1}y)\nonumber \\ &=& ... \nonumber \\
&=& (x+Q^{-(N+1)} q^{N+1} y) (x+Q^{-(N-1)} q^{N} y)...(x+Q^{N+1}y)\nonumber \\
&=& \prod_{k=0}^{N+1} (x+Q^{N+1-2k} q^k y )_{>q}. \nonumber \eea
\end{prf}
\begin{prop} For $q=1$ we have the following relation
\bea (x+y)^n = (x+y)^n_{< \tilde{Q}}, \eea
where \bea (x+y)^n_{< \tilde{Q}}= (x+Q^{-(n-1)}y) (x+Q^{-(n-3)}y)...(x+Q^{(n-3)}y) (x+Q^{(n-1)}y) \nonumber \eea is non-commutative binomial in symmetrical calculus case.
\end{prop}
\begin{prf}
By  mathematical induction, for $n=1,$ it is obvious.
Suppose we have \bea (x+y)^n = (x+y)^n_{< \tilde{Q}},\nonumber \eea for arbitrary $n.$
Then for $n+1$ we have
\bea (x+y)^{n+1} &=&(x+y)^n (x+y)=(x+y)^n_{< \tilde{Q}} (x+y)\nonumber \\
&=& (x+Q^{-(n-1)}) (x+Q^{-(n-3)})...(x+Q^{(n-3)}) (x+Q^{(n-1)}) (x+y)\nonumber \\
&=&  (x+Q^{-(n-1)}) (x+Q^{-(n-3)})...(x+Q^{(n-3)})(x+Q^{-1}y) (x+Q^{(n)})\nonumber \\
&=&  (x+Q^{-(n-1)}) (x+Q^{-(n-3)})...(x+Q^{-2}y)(x+Q^{(n-2)}) (x+Q^{(n)})\nonumber \\
&=& ... \nonumber \\
&=& (x+Q^{-n}) (x+Q^{-(n-2)})...(x+Q^{(n-2)}) (x+Q^{n})= (x+y)^{n+1}_{< \tilde{Q}}. \nonumber \eea
\end{prf}

We summarize our results in the next $q$-binomial formula for $Q$-commutative operators $x$ and $y$:
\bea (x+y)^N_{<q}&=& \prod_{k=0}^{N-1} (x+q^k y)_{<q}=(x+y) (x+ q y) (x+q^2 y)... (x+q^{N-1} y)\nonumber \\
&=& \sum_{k=0}^N {N \brack k }_{q,Q} q^{\frac{k(k-1)}{2}} x^{n-k}y^k, \label{prop1} \eea
where $yx=Qxy.$
\begin{prop} Two opposite ordered $q$-binomials are related by formula
\bea (x+y)^N_{<q}= (x+Q^{N-1}y)^N_{> \frac{q}{Q^2}},\eea where $y x= Q xy.$
\end{prop}
\begin{prf}
From  equation (\ref{relationbk}) we have
\bea \prod_{k=0}^N (x+q^k y )_{<q}=\prod_{k=0}^N (x+Q^{N-2k} q^k y )_{>q}= \prod_{k=0}^N (x+ (\frac{q}{Q^2})^k Q^N y )_{>q}= (x+Q^N y)^{N+1}_{>\frac{q}{Q^2}}; \nonumber \eea
\bea (x+y)^{N+1}_{<q}=(x+Q^N y)^{N+1}_{>\frac{q}{Q^2}}\Rightarrow \,\, (x+y)^N_{<q}= (x+Q^{N-1}y)^N_{> \frac{q}{Q^2}}.\nonumber \eea
\end{prf}

\begin{prop} For $y x= Q xy$ we obtain the following relation
\bea (x+y)^N_{>q}= \prod_{k=0}^{N-1}(x+q^k y)_{>q}=\sum_{k=0}^N {N \brack k}_{qQ^2,Q} (qQ^2)^{\frac{k(k-1)}{2}} x^{N-k} \left(\frac{y}{Q^{N-1}}\right)^k.\label{prop2}\eea
\end{prop}
\begin{prf}
We start from relation between direction of two multiplication rules (\ref{relationbk})
\bea \prod_{k=0}^{N-1} (x+q^k y )_{<q}&=&\prod_{k=0}^{N-1} (x+Q^{N-1-2k} q^k y )_{>q} \nonumber \\
&=& \prod_{k=0}^{N-1} (x+ (\frac{q}{Q^2})^k Q^{N-1} y )_{>q}.\eea
By choosing $Q^{N-1} y\equiv z\Rightarrow y= \frac{z}{Q^{N-1}},$ the above equation becomes
\bea \prod_{k=0}^{N-1}\left(x+\left(\frac{q}{Q^2}\right)^k z \right)_{>q}=\prod_{k=0}^{N-1}\left(x+q^k \frac{z}{Q^{N-1}}\right)_{<q}.\nonumber \eea
Let us call $\frac{q}{Q^2}\equiv q_1,$ then according to Proposition 2.0.10
\bea \prod_{k=0}^{N-1}(x+ q_1^k z)_{>q_1}& =& (x+z)^N_{>q_1}= \prod_{k=0}^{N-1} \left(x+(q_1 Q^2)^k \frac{z}{Q^{N-1}}\right)= \left(x+\frac{z}{Q^{N-1}}\right)^N_{<q_1 Q^2} \nonumber \\
&=& \sum_{k=0}^N {N \brack k}_{q_1 Q^2,Q} (q_1Q^2)^{\frac{k(k-1)}{2}} x^{N-k} \left(\frac{z}{Q^{N-1}}\right)^k.\eea
Relation $y x= Q xy$ implies $z x= Q xz,$ this is why, if we replace $q_1 \rightarrow q$ and $z\rightarrow y,$ we obtain the required result. \hfill \rule{1.6ex}{1.6ex}
\end{prf}
Finally, $q$-Binomial formulas for $Q$-commutative operators $x$ and $y$ with different order are summarized in equations (\ref{prop1}) and (\ref{prop2}).

Now we are going to show that all known binomial formulas like Gauss Binomial formula etc. are particular cases of our non-commutative binomial formula.
\subsection{Special Cases}

Let us consider some particular cases of this generalized $Q$ commutative $q$- Binomial formula:

(i)\,\,\, for $Q=1,$ which means commutative $x$ and $y,$ this formula becomes the Gauss Binomial formula
\bea (x+y)^n_{q} =\sum_{k=0}^n {n \brack k }_q q^{\frac{k(k-1)}{2}} x^{n-k} y^k, \nonumber \eea where $yx=xy.$

(ii)\,\,\,for $Q$- commutative $x$ and $y$ ($yx=Qxy$) and $q=1$ we have Non-commutative Binomial formula
 \bea (x+y)^n =\sum_{k=0}^n {n \brack k }_Q  x^{n-k} y^k.\nonumber \eea

(iii) \,\,\, for $Q=\frac{1}{q} ,$ we obtain the symmetrical binomial formula
\bea (x+y)^n_{q} =\sum_{k=0}^n {n \brack k }_{\tilde{q}} q^{\frac{k(k-1)}{2}} x^{n-k} y^k,\nonumber \eea where $yx=\frac{1}{q}xy,$
and the symmetrical $q$-number is defined as $$[n]_{\tilde{q}}= \frac{q^n-q^{-n}}{q-q^{-1}}.$$
(iv)\,\,\, for $q=Q \Rightarrow \lim_{Q \rightarrow q}[n]_{q,q}= n q^{n-1},$ and the formula transforms to the following one
\bea (x+y)^n_{q} =\sum_{k=0}^n {n \choose k } q^{k(n-\frac{k+1}{2})} x^{n-k} y^k,\nonumber \eea
where ${n \choose k}= \frac{n!}{(n-k)! k!}$- standard Newton binomials.

(v)\,\,\, By choosing $q=-\frac{1}{\varphi}$ and $Q=\varphi,$ where $\varphi$ is the Golden ratio, we obtain the Binet-Fibonacci Binomial formula for Golden Ratio non-commutative plane $(y x= \varphi x y)$ \cite{Pashaev and Nalci}
\bea (x+y)^n_{-\frac{1}{\varphi}} &=&\sum_{k=0}^n {n \brack k }_{\varphi,-\frac{1}{\varphi}} (-\frac{1}{\varphi})^{\frac{k(k-1)}{2}} x^{n-k} y^k \nonumber \\
&=& \sum_{k=0}^n \frac{F_n !}{F_k ! F_{n-k}!} \left(-\frac{1}{\varphi}\right)^{\frac{k(k-1)}{2}} x^{n-k} y^k ,\eea
where $F_n$ are Fibonacci numbers, and $q$-binomial coefficients become Fibonomial.

(vi) \,\,\, The above formula can be compared with the following general commutative $(Q,q)$- binomial formula. Let $x y=y x$ and 

{\small \[ (x+y)_{q,Q}^n = \left \{ \begin{array} {ll} 1 & \mbox{if $n=0$}, \\
 (x+q^{n-1} y) (x+q^{n-2}  Q y)...(x+q Q^{n-2} y)(x+Q^{n-1} y)
        & \mbox {if $n\geq 1$ } \end{array}
        \right.\] }

then the binomial formula is valid 
\bea (x+y)^n_{q,Q} =\sum_{k=0}^n {n \brack k }_{q,Q} (q Q)^{\frac{k(k-1)}{2}} x^{n-k} y^k. \nonumber \eea 
We notice that in this formula $x$ and $y$ are commutative variables, while $(q,Q)$-binomial coefficients are the same with non-commutative binomial formula.

\section{Conclusions }

In conclusion we mention one more possible application of our non-commutative $q$-binomial expansion. 
We introduce $q$-function of two variables. If 
\bea f(x)= \sum_{n=0}^\infty c_n x^n ,\nonumber \eea then $q$-function of two variables $x$ and $y$ we define as
\bea f(x+y)_q= \sum_{n=0}^\infty c_n (x+y)_q^n. \eea
These functions appear in our recent studies on $q$-analytic functions \cite{Pashaev et al.} and $q$-traveling waves \cite{cNalci et al.}.
If we take into account non-commutative binomial formulas derived in this paper we can extend our results \cite{cNalci et al.} and \cite{Pashaev et al.}  to the $q$-function of non-commutative ($Q$-commutative)  variables $x$ and $y.$ According to this we can define non-commutative $q$- analytic functions and non-commutative $q$-traveling waves. These questions are under the study.  

Here we just briefly discuss the case of non-commutative $q$-exponential function.
\begin{defn}
$(q,Q)$ analogues of exponential function are defined as
\bea e_{q,Q}(x)&\equiv& \sum_{n=0}^{\infty} \frac{1}{[n]_{q,Q}!} x^n, \nonumber \\
 E_{q,Q}(x)&\equiv & \sum_{n=0}^{\infty} \frac{1}{[n]_{q,Q}!}q^{\frac{n(n-1)}{2}} x^n.\eea
\end{defn}
\begin{prop} For $Q$-commutative operators $x$ and $y,$ $(yx=Qxy),$ we have the following factorization of $q$- exponential function $e_{q,Q},$
$$e_{q,Q} (x+y)_{<q}= e_{q,Q} (x) E_{q,Q} (y),$$
\vspace{-1cm}
\end{prop} 
\vspace{+1cm}
\begin{prf}
\bea e_{q,Q} (x+y)_{<q}&=& \sum_{N=0}^{\infty} \frac{(x+y)^N_{<q}}{[N]_{q,Q}!}\nonumber \\
&=& \sum_{N=0}^{\infty}\frac{1}{[N]_{q,Q}!} \sum_{k=0}^N {N \brack k}_{q,Q} q^{\frac{k(k-1)}{2}} x^{N-k} y^k \nonumber \\
&=&  \sum_{N=0}^{\infty}\sum_{k=0}^N \frac{1}{[N-k]_{q,Q}![k]_{q,Q}!} q^{\frac{k(k-1)}{2}} x^{N-k} y^k. \nonumber \eea
By choosing $N-k \equiv s,$
\bea e_{q,Q} (x+y)_{<q}&=& \left( \sum_{s=0}^{\infty} \frac{1}{[s]_{q,Q}!} x^s \right) \left( \sum_{k=0}^{\infty} \frac{1}{[k]_{q,Q}!}q^{\frac{k(k-1)}{2}} y^k \right)\nonumber \\
&=& e_{q,Q} (x) E_{q,Q} (y).\eea \hfill \rule{1.6ex}{1.6ex}

\end{prf}

For $q=1,$ the above exponential functions reduce to the Jackson exponential functions and our proposition gives factorization of this function for $Q$-commutative argument.

 \section*{Acknowledgments}
This work was carried out with support from
TUBITAK (The Scientific and Technological Research Council of Turkey), TBAG Project 110T679 and Izmir Institute of
Technology. And one of the authors (S.Nalci) was partially supported by TUBITAK scholarship of graduate students.

\end{document}